# The upper bound of the spectral radius for the hypergraphs without Berge-graphs

Wen-Huan Wang[*], Lou-Jun Yu

*Department of Mathematics, Shanghai University, Shanghai 200444, China*

(Dated: September 27, 2023)

The spectral analogue of the Turán type problem for hypergraphs is to determine the maximum spectral radius for the hypergraphs of order $n$ that do not contain a given hypergraph. For the hypergraphs among the set of the connected linear 3-uniform hypergraphs on $n$ vertices without the Berge–$C_l$, we present two upper bounds for their spectral radius and $\alpha$-spectral radius, which are related to $n$, $l$ and $\alpha$, where $C_l$ is a cycle of length $l$ with $l \geqslant 5$, $n \geqslant 3$ and $0 \leqslant \alpha < 1$. Let $B_s$ be an $s$-book with $s \geqslant 2$ and $K_{s,t}$ be a complete bipartite graph with two parts of size $s$ and $t$, respectively, where $s, t \geqslant 1$. For the hypergraphs among the set of the connected linear $k$-uniform hypergraphs on $n$ vertices without the Berge–$\{B_s, K_{2,t}\}$, we derive two upper bounds for their spectral radius and $\alpha$-spectral radius, which depend on $n$, $k$, $s$, and $\alpha$, where $n, k \geqslant 3$, $s \geqslant 2$, $1 \leqslant t \leqslant \frac{1}{2}(6k^2 - 15k + 10)(s - 1) + 1$, and $0 \leqslant \alpha < 1$.

## 1. INTRODUCTION

Let $\mathbb{R}$ and $\mathbb{C}$ be the sets of real and complex numbers, respectively. A $k$-ordered and $n$-dimensional real tensor $\boldsymbol{\mathcal{A}} = (a_{i_1 i_2 \cdots i_k})$ over $\mathbb{R}$ is a multi-dimensional array with $n^k$ entries, where $a_{i_1 i_2 \cdots i_k} \in \mathbb{R}$ with $i_1, i_2, \cdots, i_k \in [n] = \{1, 2, \cdots, n\}$. The concept of tensor eigenvalues and the spectra of tensors are independently introduced by Qi [1] and Lim [2] as follows. If there exist a number $\lambda \in \mathbb{C}$ and an eigenvector $\boldsymbol{x} = \{x_1, x_2, \cdots, x_n\}^{\mathrm{T}} \in \mathbb{C}^n$

---

[*] Corresponding author. Email: whwang@shu.edu.cn

satisfying

$$\sum_{i_2,\cdots,i_k=1}^{n} a_{ii_2\cdots i_k}x_{i_2}\cdots x_{i_k} = \lambda x_i^{k-1}, \text{ for any } 1 \leqslant i \leqslant n, \tag{1}$$

then $\lambda$ is called an eigenvalue of $\boldsymbol{\mathcal{A}}$ and $\boldsymbol{x}$ an eigenvector of $\boldsymbol{\mathcal{A}}$ corresponding to $\lambda$. The spectral radius of $\boldsymbol{\mathcal{A}}$, denoted by $\rho(\boldsymbol{\mathcal{A}})$, is the largest modulus of all the eigenvalues of $\boldsymbol{\mathcal{A}}$.

A hypergraph $H$ is a pair $(V(H), E(H))$, where $V(H)$ and $E(H)$ stand for the sets of the vertices and the edges of $H$, respectively. For an arbitrary edge $e \in E(H)$, if $|e| = k$, then $H$ is a $k$-uniform hypergraph, where $k \geqslant 2$. For simplicity, a $k$-uniform hypergraph is hereinafter called a $k$-graph, where $k \geqslant 2$. Obviously, a 2-graph is a graph. A $k$-graph $H$ is linear if any two edges of $E(H)$ intersect at most one vertex. For $u, v \in V(H)$, if there exists an edge $e \in E(H)$ such that $\{u, v\} \subseteq e$, then $u$ and $v$ are adjacent, and $u$ is incident with $e$. The number of the edges of $H$ incident with $v$, denoted by $d_H(v)$, is the degree of $v$. Without confusion, $d_H(v)$ is simplified as $d_v$.

Let $H = (V(H), E(H))$ be a $k$-graph on $n$ vertices, where $V(H) = \{v_1, \cdots, v_n\}$ and $k, n \geqslant 2$. The adjacency tensor of $H$ is $\boldsymbol{\mathcal{A}}(H) = (a_{i_1i_2\cdots i_k})$, where $a_{i_1i_2\cdots i_k} = \frac{1}{(k-1)!}$ if $\{v_{i_1}, v_{i_2}, \cdots, v_{i_k}\} \in E(H)$ and $a_{i_1i_2\cdots i_k} = 0$ otherwise [3]. The degree diagonal tensor of $H$ is $\boldsymbol{\mathcal{D}}(H) = (d_{i_1i_2\cdots i_k})$, where $d_{i_1i_2\cdots i_k} = d_{v_i}$ for any $v_i \in V(H)$ if $i_1 = i_2 = \cdots = i_k = i$ with $i \in [n]$ and $d_{i_1i_2\cdots i_k} = 0$ otherwise. Let $0 \leqslant \alpha < 1$. Lin et al. [4] introduced the convex linear combination of $\boldsymbol{\mathcal{D}}(H)$ and $\boldsymbol{\mathcal{A}}(H)$ for a hypergraph $H$ as follows: $\boldsymbol{\mathcal{A}}_\alpha(H) = \alpha\boldsymbol{\mathcal{D}}(H) + (1-\alpha)\boldsymbol{\mathcal{A}}(H)$. The $\alpha$-spectral radius of $H$, denoted by $\rho_\alpha(H)$, is the spectral radius of $\boldsymbol{\mathcal{A}}_\alpha(H)$. Obviously, when $\alpha = 0$, $\boldsymbol{\mathcal{A}}_\alpha(H)$ is the adjacency tensor and $\rho_0(H)$ is equal to the spectral radius of $H$ (denoted by $\rho(H)$). When $\alpha = \frac{1}{2}$, $2\boldsymbol{\mathcal{A}}_\alpha(H)$ is the signless Laplacian tensor of $H$ and $2\rho_{\frac{1}{2}}(H)$ is the signless Laplacian spectral radius of $H$.

Let $\mathcal{F}$ be a given family of hypergraphs. If a hypergraph $H$ does not contain any member of $\mathcal{F}$ as a subhypergraph, then $H$ is $\mathcal{F}$–free. If $\mathcal{F} = \{F\}$, then $H$ is $F$–free. The Turán type extremal problem is one of the most fundamental problems in the extremal combinatorics, which is to determine the Turán number of $\mathcal{F}$, namely the maximum number of edges in $\mathcal{F}$–free hypergraphs on $n$ vertices. For the results about the Turán type extremal problems for graphs and hypergraphs, one can refer to surverys [5, 6] and Refs. [7–12], respectively.

Let $G = (V(G), E(G))$ be a given graph and $H = (V(H), E(H))$ be a hypergraph. $H$ is called a Berge–$G$ if there exists a bijection $\varphi : E(G) \to E(H)$ such that $e \subseteq \varphi(e)$ for any edge $e \in E(G)$ [9]. Namely, $H$ is obtained from $G$ by replacing each edge of $G$ with a hyperedge that contains it. For a given graph $G$, a hypergraph $H$ is Berge–$G$–free if $H$ does not contain a subhypergraph isomorphic to a Berge–$G$. Recently, the study on the Turán numbers of hypergraphs without a Berge–$G$ attracted the attention of researchers [13, 14].

The spectral version of the Turán type problem for graphs is to determine the maximum spectral radius among all the $G$–free graphs on $n$ vertices, which was first proposed by Nikiforov [15]. For details, one can refer to Refs. [16–21].

For the spectral analogue of the Turán type problem for hypergraphs, Hou et al. [22] proposed the following question: "What is the maximum spectral radius of hypergraphs of order $n$, not containing a given hypergraph $F$?" This problem attracts researchers' interest and some results are obtained. eevash et al. [23] determined the maximum $p$-spectral radius of any 3-graphs on $n$ vertices not containing the Fano plane when $n$ is sufficently large. Hou et al. [22] derived the upper bounds of the maximum spectral radii for the connected linear $k$-graphs on $n$ vertices without the Berge–$C_4$ or with the girth of at least five, respectively. For a given graph $G$, the $r$-expansion $G^+$ of $G$ is the $r$-uniform hypergraph obtained from $G$ by enlarging each edge of $G$ with $r-2$ new vertices which are disjoint from the vertex set of $G$ such that distinct edges are enlarged by disjoint subsets. Gao et al. [24] obtained the maximum spectral radius of the connected linear $r$-graphs on $n$ vertices without $K_{r+1}^+$, where $K_{r+1}$ is a complete graph on $r + 1$ vertices. She and Fan [25] presented the sharp (or asymptotic) bounds of the maximum number of edges and the maximum spectral radius of all $F^+$–free linear $r$-uniform hypergraphs on $n$ vertices, respectively, where $F$ is a graph with chromatic number $k$.

Let $P_k$ be a path of length $k$ and we write $P_k$ as $v_1 v_2 \cdots v_{k+1}$, where $k \geqslant 1$. A cycle of length $k$, denoted by $C_k$, is obtained from $P_k$ by identifying $v_1$ with $v_{k+1}$ together. Let $sC_3$ be the union of $C_3$ and the number of $C_3$ is $s$, where $s \geqslant 2$. An $s$-book, denoted by $B_s$, is a connected graph obtained from $sC_3$ in such a way that all $C_3$ share a common edge. Let $K_{s,t}$ be a complete bipartite graph with two parts of size $s$ and $t$, respectively, where $s, t \geqslant 1$.

The research on the spectral analogue of the Turán-type problem for hypergraphs is an interesting field. Let $\mathcal{H}_{n,3,l}$ be the set of the connected linear 3-graphs on $n$ vertices without the Berge–$C_l$, where $n \geqslant 3$ and $l \geqslant 5$. Let $\mathcal{H}_{n,k,s,t}$ be the set of the connected linear $k$-graphs on $n$ vertices without the Berge–$\{B_s, K_{2,t}\}$, where $n, k \geqslant 3$, $s \geqslant 2$ and $t \geqslant 1$. In this article, for the hypergraphs among $\mathcal{H}_{n,3,l}$ and $\mathcal{H}_{n,k,s,t}$, we aim to characterize the upper bounds of their maximum spectral radius and $\alpha$-spectral radius.

The rest of this article is organized as follows. In Section 2, some known necessary lemmas are quoted and some new tools for us to study the $\alpha$-spectral radius of hypergraphs are developed. In Section 3, for $H \in \mathcal{H}_{n,3,l}$ with $n \geqslant 3$ and $l \geqslant 5$, we obtain the upper bounds for $\rho(H)$ and $\rho_\alpha(H)$, which depend on $n$, $l$ and $\alpha$. In Section 4, for $H \in \mathcal{H}_{n,k,s,t}$ with $n, k \geqslant 3$, $s \geqslant 2$ and $1 \leqslant t \leqslant \frac{1}{2}(6k^2 - 15k + 10)(s-1) + 1$, the upper bounds for $\rho(H)$ and $\rho_\alpha(H)$ are also derived, which are related with $n$, $k$, $s$, and $\alpha$.

## 2. PRELIMINARIES

In this section, we will introduce several necessary lemmas for subsequent proofs.

The nonnegative weakly irreducible tensor was defined by Friedland et al. [26] and was represented by Yang et al. [27] as follows.

**Definition 2.1** *[26, 27] Let $\mathcal{A} = (a_{i_1 i_2 \cdots i_k})$ be a nonnegative tensor of order $k$ and dimension $n$. For any nonempty proper index subset $I \subset [n]$, if there is at least one entry $a_{i_1 i_2 \cdots i_k} > 0$, where $i_1 \in I$ and at least one $i_j \in [n] \setminus I$ for $j = 2, 3, \cdots, k$, then $\mathcal{A}$ is called a nonnegative weakly irreducible tensor.*

Let $\mathbb{R}_+^n = \{\boldsymbol{x} \in \mathbb{R}^n \mid x_i \geqslant 0, \forall i \in [n]\}$ and $\mathbb{R}_{++}^n = \{\boldsymbol{x} \in \mathbb{R}^n \mid x_i > 0, \forall i \in [n]\}$.

**Lemma 2.1** *[26, 28] (The Perron–Frobenius theorem for nonnegative tensors) Let $\mathcal{A}$ be a nonnegative tensor of order $k$ and dimension $n$, where $k, n \geqslant 2$. Then we have the following statements.*

*(i). $\rho(\mathcal{A})$ is an eigenvalue of $\mathcal{A}$ with a nonnegative eigenvector $\boldsymbol{x} \in \mathbb{R}_+^n$ corresponding to it.*

*(ii). If $\mathcal{A}$ is weakly irreducible, then $\rho(\mathcal{A})$ is the unique eigenvalue of $\mathcal{A}$ with the positive eigenvector $\boldsymbol{x} \in \mathbb{R}_{++}^n$, up to a positive scaling coefficient.*

**Lemma 2.2** *[29] A $k$-graph $H$ is connected if and only if $\mathcal{A}(H)$ is weakly irreducible.*

If $\mathcal{A}(H)$ is weakly irreducible, then according to the definition of $\mathcal{A}_\alpha(H)$, $\mathcal{A}_\alpha(H)$ is also weakly irreducible. If $H$ is a connected $k$-graph, then by Lemmas 2.1 and 2.2, there exists a unique vector $\boldsymbol{x} = (x_1, x_2, \cdots, x_n)^{\mathrm{T}}$ corresponding to $\rho_\alpha(H)$ and each $x_i > 0$ with $1 \leqslant i \leqslant n$. This vector $\boldsymbol{x}$ is referred to as the $\alpha$-Perron vector of $H$. It is convenient for us to normalize $\boldsymbol{x}$ such that the maximum entry of $\boldsymbol{x}$ is 1. If there are many such vertices, then among them we can arbitrarily choose one and fix it. For any $v \in V(H)$, from (1), we have

$$\rho_\alpha(H) x_v^{k-1} = \alpha d_v x_v^{k-1} + (1-\alpha) \sum_{\{v, v_{i_2}, \cdots, v_{i_k}\} \in E(H)} x_{v_{i_2}} x_{v_{i_3}} \cdots x_{v_{i_k}}. \tag{2}$$

**Lemma 2.3** *[30] (Erdős–Gallai Theorem) Let $G$ be a graph on $n$ vertices and $e(G)$ be the number of edges of $G$. If $G$ is $P_k$–free, then $e(G) \leqslant \frac{k-1}{2} n$ with the equality if and only if $G$ is composed of multiple $K_k$ which are disjoint from each other.*

Let $H$ be a $k$-graph and let $v \in V(H)$, where $k \geqslant 2$. For a subset $X \subseteq V(H)$, let $E_t(X) = \{e \mid e \in E(H) \text{ and } |e \cap X| = t\}$, where $t \geqslant 1$. Namely, each edge in $E_t(X)$ and $X$ intersect at $t$ vertices. We define $|E_t(X)| = e_t(X)$. Furthermore, let $E_t^v(X) = \{e \mid e \in E_t(X) \text{ and } v \in e\}$ and $|E_t^v(X)| = e_t^v(X)$. Namely, $E_t^v(X)$ is a subset of $E_t(X)$ satisfying that each edge in $E_t^v(X)$ contains $v$. We define $N_v = \{w \mid w \in V(H), \{v, w\}$ is contained in an edge of $H\}$ and call $N_v$ the neighbourhood of $v$. In $H$, let $V(H) = \{v_1, v_2, \cdots, v_n\}$. A Berge walk of length $h$ in $H$ is an alternating sequence of vertices and edges $v_1 e_1 v_2 e_2 v_3 \cdots v_h e_h v_{h+1}$, where $v_i \neq v_{i+1}$ and $\{v_i, v_{i+1}\} \subseteq e_i$ for $1 \leqslant i \leqslant h$. We use $w_h(v)$ to denote the number of the Berge walks of length $h$ in $H$ which start at a vertex $v \in V(H)$. By the definition of $w_h(v)$, we have

$$w_1(v) = (k-1)d_v, \tag{3}$$

$$w_{h+1}(v) = \sum_{\{v, v_{i_2}, \cdots, v_{i_k}\} \in E(H)} \left[ w_h(v_{i_2}) + w_h(v_{i_3}) + \cdots + w_h(v_{i_k}) \right]. \tag{4}$$

**Lemma 2.4** *[24] Let $H$ be a connected linear $k$-graph, where $k \geqslant 3$. For any $v \in V(H)$, we have*

$$w_2(v) = (k-1) \sum_{i=1}^{k} i e_i(N_v).$$

To obtain the upper bound of the $\alpha$-spectral radius for the hypergraphs, we introduce some new tools, which are shown in Theorem 2.1, Corollaries 2.1 and 2.2.

**Theorem 2.1** *Let $H = (V(H), E(H))$ be a connected $k$-graph on $n$ vertices and $\boldsymbol{x} = (x_1, \cdots, x_n)^{\mathrm{T}}$ be its $\alpha$-Perron vector with the maximum entry in $\boldsymbol{x}$ being 1, where $k, n \geqslant 3$. For any $v \in V(H)$, we have*

$$\rho_\alpha^2(H) x_v^{k-1} \leqslant \frac{\alpha}{(k-1)^2} w_1^2(v) + \frac{1-\alpha}{(k-1)^2} w_2(v).$$

**Proof.** Let $v$ be an arbitrary vertex in $V(H)$. Since $x_i \leqslant 1$ for any $1 \leqslant i \leqslant n$, by (2), we get

$$
\begin{aligned}
\rho_\alpha(H)x_v^{k-1} &= \alpha d_v x_v^{k-1} + (1-\alpha) \sum_{\{v,i_2,\cdots,i_k\}\in E(H)} x_{i_2}\cdots x_{i_k} \\
&\leqslant \alpha d_v + (1-\alpha) \sum_{\{v,i_2,\cdots,i_k\}\in E(H)} 1 \\
&= d_v = \frac{1}{k-1}w_1(v),
\end{aligned} \tag{5}
$$

where (5) follows from (3). We obtain

$$
\begin{aligned}
\rho_\alpha^2(H)x_v^{k-1} &= \rho_\alpha(H)\cdot\rho_\alpha(H)x_v^{k-1} \\
&= \rho_\alpha(H)\left[\alpha d_v x_v^{k-1} + (1-\alpha) \sum_{\{v,i_2,\cdots,i_k\}\in E(H)} x_{i_2}\cdots x_{i_k}\right] \\
&\leqslant \alpha d_v \rho_\alpha(H)x_v^{k-1} + \frac{1-\alpha}{k-1}\sum_{\{v,i_2,\cdots,i_k\}\in E(H)}\rho_\alpha(H)(x_{i_2}^{k-1}+\cdots+x_{i_k}^{k-1}) \tag{6} \\
&\leqslant \alpha d_v^2 + \frac{1-\alpha}{(k-1)^2}\sum_{\{v,i_2,\cdots,i_k\}\in E(H)}\left(w_1(i_2)+\cdots+w_1(i_k)\right) \tag{7} \\
&= \frac{\alpha}{(k-1)^2}w_1^2(v) + \frac{1-\alpha}{(k-1)^2}w_2(v), \tag{8}
\end{aligned}
$$

where (6) follows from the Arithmetic Mean–Geometry Mean inequality, (7) is derived from (5), and (8) is deduced from (3) and (4). $\qquad\square$

By Theorem 2.1, we directly get Corollary 2.1.

**Corollary 2.1** *Let $H = (V(H), E(H))$ be a connected $k$-graph on $n$ vertices and $\boldsymbol{x} = (x_1,\cdots,x_n)^{\mathrm{T}}$ be its $\alpha$-Perron vector, where $k, n \geqslant 3$. If $u \in V(H)$ is the vertex having the maximum entry in $\boldsymbol{x}$, then*

$$
\rho_\alpha^2(H) \leqslant \frac{\alpha}{(k-1)^2}w_1^2(u) + \frac{1-\alpha}{(k-1)^2}w_2(u).
$$

**Corollary 2.2** *Let $H = (V(H), E(H))$ be a connected linear $k$-graph on $n$ vertices and $\boldsymbol{x} = (x_1,\cdots,x_n)^{\mathrm{T}}$ be its $\alpha$-Perron vector, where $k, n \geqslant 3$. If $u \in V(H)$ is the vertex having the maximum entry in $\boldsymbol{x}$, then*

$$
\rho_\alpha^2(H) \leqslant \alpha d_u^2 + \frac{1-\alpha}{k-1}\sum_{i=1}^k ie_i(N_u) = \alpha d_u^2 + \frac{1-\alpha}{k-1}\sum_{v\in N_u} d_v. \tag{9}
$$

**Proof.** By combining Corollary 2.1, Lemma 2.4, (3), and $\sum_{i=1}^k ie_i(N_u) = \sum_{v\in N_u} d_v$, we obtain Corollary 2.2. $\qquad\square$

It is noted that (9) with $\alpha = 0$ was obtained by Hou et al. [22] and Gao et al. [24]. It can be seen from Corollary 2.2 that we can use a unified method to study the upper bounds of the maximum spectral radius and $\alpha$-spectral radius for a connected linear $k$-graph on $n$ vertices.

## 3. THE UPPER BOUND OF THE MAXIMUM SPECTRAL RADIUS OF THE CONNECTED LINEAR 3-GRAPHS WITHOUT THE BERGE–$C_l$

In this section, we will characterize the upper bounds of the maximumspectral radius and $\alpha$-spectral radius for the hypergraphs among the set of the connected linear 3-graphs on $n$ vertices without the Berge–$C_l$, where $n \geqslant 3$ and $l \geqslant 5$.

If $H$ is a linear hypergraph, for any two adjacent vertices $x$ and $y$ of $H$, there exists a unique edge containing $x$ and $y$. We denote this edge by $l_{xy}$.

**Theorem 3.1** Let $H \in \mathcal{H}_{n,3,l}$, where $n \geqslant 3$ and $l \geqslant 5$. Let $\boldsymbol{x} = (x_1, \cdots, x_n)^{\mathrm{T}}$ be the $\alpha$-Perron vector of $H$. If $u \in V(H)$ is the vertex having the maximum entry in $\boldsymbol{x}$, then

$$\rho_\alpha^2(H) \leqslant \alpha d_u^2 + (1 - \alpha) \left[ \frac{(2l - 7)(n - 1)}{4} + (2l - 3)d_u \right].$$

**Proof.** To obtain Theorem 3.1, we firstly prove Claims 3.1 and 3.2.

Let $B_u = V(H) \backslash (N_u \cup \{u\})$. Since $|N_u| = 2d_u$, we get $|B_u| = n - 2d_u - 1$.

**Claim 3.1.** $e_1(N_u) \leqslant \frac{2l-7}{2}(n - 2d_u - 1)$.

**Proof.** We will prove Claim 3.1 by constructing an auxiliary graph. Let $G_1 = \big(V(G_1), E(G_1)\big)$, where $V(G_1) = B_u$ and $E(G_1) = \{v_1 v_2 \mid v_1, v_2 \in B_u, \ l_{v_1 v_2} \in E_1(N_u)\}$. By the definition of $G_1$, we have $|V(G_1)| = |B_u| = n - 2d_u - 1$ and $|E(G_1)| = e_1(N_u)$.

We prove $|E(G_1)| \leqslant \frac{2l-7}{2}(n - 2d_u - 1)$ by contradiction. Otherwise, we suppose $|E(G_1)| > \frac{2l-7}{2}(n - 2d_u - 1) = \frac{2l-7}{2}|V(G_1)|$. Obviously, the average degree of $G_1$ is greater than $2l - 7$. When $l \geqslant 5$, we have $2l - 7 \geqslant 3$. In $G_1$, if we delete the vertex having degree less than 3, then the average degree of the resulting graph does not decrease. Therefore, in $G_1$, if we delete all the vertices having degree less than 3, then we get a resulting graph (denoted by $G_2$) such that the average degree of $G_2$ is greater than $2l - 7$ and the minimum degree of $G_2$ is greater than or equal to 3. For $G_2$, since its average degree is greater than $2l - 7$, we get $|E(G_2)| > \frac{2l-7}{2}|V(G_2)|$. By Lemma 2.3, we obtain that $G_2$ contains a path $P_{2l-6}$ as a subgraph, where $l \geqslant 5$. Let $P_{2l-6} = v_1 e_1 v_2 \cdots v_{2l-6} e_{2l-6} v_{2l-5}$.

Since the minimum degree of $G_2$ is greater than or equal to 3, there exists a vertex (denoted by $w$) in $B_u$ such that $wv_{l-2} \in E(G_2)$, where $w \neq v_{l-3}, v_{l-1}$. According to the definition of $G_1$, in $H$, there exist five vertices $w_1, \cdots, w_5 \in N_u$ and five edges $f_1, \cdots, f_5$, where $f_1 = \{v_1, v_2, w_1\}$, $f_2 = \{v_{l-3}, v_{l-2}, w_2\}$, $f_3 = \{v_{l-2}, v_{l-1}, w_3\}$, $f_4 = \{v_{2l-6}, v_{2l-5}, w_4\}$, and $f_5 = \{v_{l-2}, w, w_5\}$. It is noted that $w_1, w_2, \cdots, w_5$ may be the same. Thus, we obtain $l_{uw_1} = l_{uw_3}$. Otherwise, if $l_{uw_1} \neq l_{uw_3}$, then there exists a Berge–$C_l$ in $H$ which is formed by $l_{uw_1}$, $l_{w_1v_2}$, $l_{v_2v_3}$, $\cdots$, $l_{v_{l-2}w_3}$, $l_{uw_3}$. This contradicts the fact that $H$ is Berge–$C_l$–free. By the same analysis as those for $l_{uw_1} = l_{uw_3}$, we also get $l_{uw_1} = l_{uw_5}$. Therefore, we have $l_{uw_1} = l_{uw_3} = l_{uw_5}$. Furthermore, we claim $l_{uw_2} = l_{uw_4}$. Otherwise, if $l_{uw_2} \neq l_{uw_4}$, then there exists a Berge–$C_l$ in $H$ which is formed by $l_{uw_2}$, $l_{w_2v_{l-2}}$, $l_{v_{l-2}v_{l-1}}$, $\cdots$, $l_{v_{2l-6}w_4}$, $l_{uw_4}$. This contradicts the fact that $H$ is Berge–$C_l$–free. Similarly, we get $l_{uw_4} = l_{uw_5}$. Therefore, we obtain $l_{uw_2} = l_{uw_4} = l_{uw_5}$. Thus, $l_{uw_2} = l_{uw_3} = l_{uw_5}$. This means that an edge of $H$ contains $u$, $w_2$, $w_3$, and $w_5$. Since $v_{l-2} \in f_2 \cap f_3 \cap f_5$, $w_2 \in f_2$, $w_3 \in f_3$, $w_5 \in f_5$, and $H$ is a linear hypergraph, we deduce that $w_2$, $w_3$, and $w_5$ are different from each other. Therefore, an edge of $H$ has at least four distinct vertices. This contradicts the fact that $H$ is a 3-graph. Thus, we get $e_1(N_u) = |E(G_1)| \leqslant \frac{2l-7}{2}(n - 2d_u - 1)$. Namely, Claim 3.1 holds. $\qquad \square$

**Claim 3.2.** $2e_2(N_u) + 3e_3(N_u) \leqslant (6l - 13)d_u$.

**Proof.** Since $H$ is a 3-graph, if an edge of $H$ is incident with $u$, then this edge and $N_u$ intersect at two vertices. In $H$, the number of such edges which, containing $u$, intersect with $N_u$ at two vertices is $d_u$. Namely, $|E_2^u(N_u)| = d_u$.

Next, we will prove the following Fact 1 by constructing an auxiliary graph.

In $H$, let $E'$ be the set of such edges which, not containing $u$, intersect with $N_u$ at two vertices or three vertices.

**Fact 1**: $|E'| \leq (2l - 5)d_u$.

Let $G_3 = \left(V(G_3), E(G_3)\right)$, where $V(G_3) = N_u$ and the edges of $E(G_3)$ are formed by the edges of $E'$ according to the following one-to-one mapping relationship. Let $e$ be an arbitrary edge of $E'$. If $|e \cap N_u| = 2$, namely, there exist two vertices $x, y \in N_u$ such that $e \cap N_u = \{x, y\}$, then let $xy \in E(G_3)$. If $|e \cap N_u| = 3$, namely, there exist three vertices $x, y, z \in N_u$ such that $e \cap N_u = \{x, y, z\}$, then we take any two vertices from $\{x, y, z\}$ to form an edge of $G_3$. Obviously, according to the mapping relationship, each edge in $E'$ corresponds to a unique edge of $E(G_3)$. Therefore, we get $|E(G_3)| = |E'|$.

Next, we prove that $G_3$ does not contain $P_{2l-4}$. Otherwise, if $G_3$ contains $P_{2l-4} = v_1e_1v_2\cdots v_{2l-4}e_{2l-4}v_{2l-3}$, then by the construction of $G_3$, for any $i = 1, 2, \cdots, 2l-3$, we get $v_i \in N_u$. Since $v_{l-1} \in N_u$, in $H$, there exists an edge $f_1 \in E(H)$ such that $\{u, v_{l-1}\} \subseteq f_1$. As $H$ is a 3-graph, $v_1 \in f_1$ and $v_{2l-3} \in f_1$ can not hold simultaneously. Without loss of generality, let $v_1 \notin f_1$. Since $v_1 \in N_u$, in $H$, there exists an edge $f_2$ such that $\{u, v_1\} \subseteq f_2$ and $f_2 \neq f_1$. Therefore, in $H$, $f_1$, $f_2$, $l_{v_1v_2}, \cdots, l_{v_{l-2}v_{l-1}}$ form a Berge–$C_l$. This contradicts the fact that $H$ is Berge–$C_l$–free. Therefore, $G_3$ is $P_{2l-4}$–free. By Lemma 2.3, we get $|E(G_3)| \leqslant \frac{2l-4-1}{2}|V(G_3)| = \frac{2l-5}{2}|N_u|$. Since $H$ is a connected linear 3-graph, $|N_u| = 2d_u$. Thus, $|E'| = |E(G_3)| \leqslant (2l-5)d_u$. Namely, Fact 1 holds.

Bearing $|E_2^u(N_u)| = d_u$ and Fact 1 in mind, we obtain

$$
\begin{aligned}
2e_2(N_u) + 3e_3(N_u) &= 2|E_2(N_u)| + 3|E_3(N_u)| \\
&= 2|E_2^u(N_u)| + 2|E_2(N_u)\backslash E_2^u(N_u)| + 3|E_3(N_u)| \\
&\leqslant 2|E_2^u(N_u)| + 3\Big(|E_2(N_u)\backslash E_2^u(N_u)| + |E_3(N_u)|\Big) \\
&= 2|E_2^u(N_u)| + 3|E'| \leqslant (6l-13)d_u.
\end{aligned}
\tag{10}
$$

Since $H$ is a connected linear 3-graph, by Corollary 2.2, Claims 3.1 and 3.2, and (10), we have

$$
\begin{aligned}
\rho_\alpha^2(H) &\leqslant \alpha d_u^2 + \frac{1-\alpha}{2} \sum_{i=1}^{3} ie_i(N_u) \\
&\leqslant \alpha d_u^2 + \frac{1-\alpha}{2} \left[\frac{2l-7}{2}(n - 2d_u - 1) + (6l-13)d_u\right] \\
&= \alpha d_u^2 + (1-\alpha)\left[\frac{(2l-7)(n-1)}{4} + (2l-3)d_u\right].
\end{aligned}
\tag{11}
$$

Thus, we get Theorem 3.1. □

We can see from (11) that the upper bound of $\rho_\alpha(H)$ has a close relationship with $d_u$ when $n$, $l$ and $\alpha$ are fixed. Therefore, by Theorem 3.1, we directly obtain Theorem 3.2.

**Theorem 3.2** *Let $H$ be a connected linear 3-graph on $n \geqslant 3$ vertices with the maximum degree $\Delta$. If $H$ is Berge–$C_l$–free, where $l \geqslant 5$, then*

$$
\rho_\alpha^2(H) \leqslant \alpha\Delta^2 + (1-\alpha)\left[\frac{(2l-7)(n-1)}{4} + (2l-3)\Delta\right].
$$

If $H$ is a connected linear 3-graph on $n \geqslant 3$ vertices, then its maximum degree $\Delta$ is less than or equal to $\frac{n-1}{2}$. Thus, by Theorem 3.2, for $H \in \mathcal{H}_{n,3,l}$ with $n \geqslant 3$ and $l \geqslant 5$, we obtain the upper bounds of $\rho(H)$ and $\rho_\alpha(H)$, which are shown in Theorems 3.3 and 3.4, respectively.

**Theorem 3.3** *Let $H$ be a connected linear 3-graph on $n$ vertices with $n \geqslant 3$. If $H$ is Berge–$C_l$–free, where $l \geqslant 5$, then*

$$\rho(H) \leqslant \sqrt{\frac{(6l-13)(n-1)}{4}}.$$

**Theorem 3.4** *Let $H$ be a connected linear 3-graph on $n$ vertices with $n \geqslant 3$. If $H$ is Berge–$C_l$–free, where $l \geqslant 5$, then*

$$\rho_\alpha^2(H) \leqslant \alpha \frac{(n-1)^2}{4} + (1-\alpha)\frac{(6l-13)(n-1)}{4}.$$

# 4. THE UPPER BOUND OF THE MAXIMUM SPECTRAL RADIUS FOR THE $k$-UNIFORM LINEAR HYPERGRAPHS WITHOUT THE BERGE–$\{B_s, K_{2,t}\}$

In this section, we will derive the upper bounds of the maximum spectral radius and $\alpha$-spectral radius for the hypergraphs among the set of the connected linear $k$-graphs on $n$ vertices without the Berge–$\{B_s, K_{2,t}\}$, where $n, k \geqslant 3$, $s \geqslant 2$ and $t \geqslant 1$.

**Theorem 4.1** *Let $H \in \mathcal{H}_{n,k,s,t}$, where $k, n \geqslant 3$, $s \geqslant 2$, and $t \geqslant 1$. Let $\boldsymbol{x} = (x_1, \cdots, x_n)^{\mathrm{T}}$ be the $\alpha$-Perron vector of $H$. If $u \in V(H)$ is the vertex having the maximum entry in $\boldsymbol{x}$, then*

$$\rho_\alpha^2(H) \leqslant \alpha d_u^2 + (1-\alpha)\left[\frac{(t-1)(n-1)}{(k-1)^2} + \frac{(6k^2 - 15k + 10)(s-1) - 2(t-1)}{2(k-1)}d_u\right]. \tag{12}$$

**Proof.** To obtain Theorem 4.1, we firstly prove Claims 4.1–4.4.

For two vertices $x$ and $y$ in $H$, let $N_{xy}$ be the set of the common neighbourhood of $x$ and $y$. Namely, each vertex in $N_{xy}$ is adjacent to $x$ and $y$ simultaneously. Furthermore, we define $N'_{xy} = N_{xy} \backslash l_{xy}$.

**Claim 4.1.** For any two adjacent vertices $x, y \in V(H)$, we have $|N_{xy}| \leqslant (k-1)(2s-1) - s$.

**Proof.** We prove Claim 4.1 by contradiction. Suppose that there exist two adjacent vertices $x, y \in V(H)$ such that $|N_{xy}| \geqslant (k-1)(2s-1) - s + 1 = (2k-3)(s-1) + k - 1$.

Then $|N'_{xy}| \geqslant (2k-3)(s-1)+1$. Let $V_0 = N'_{xy}$ and $E_0 = \{l_{xy}\}$. Next, we will prove that there exists a Berge-$B_s$ in $H$ according to the following recursive method.

Step (1): Since $V_0 = N'_{xy} \neq \phi$, we can choose one vertex from $V_0$ and denote this vertex by $w_1$. Let $E_1 = E_0 \cup \{l_{xw_1}, l_{yw_1}\}$ and $V_1 = V_0 \backslash (l_{xw_1} \cup l_{yw_1})$. Since $x, y \notin V_0$ and $|(l_{xw_1} \cup l_{yw_1}) \cap V_0| \leqslant 2k-3$, we get $|V_1| \geqslant |V_0| - (2k-3) = |N'_{xy}| - (2k-3) \geqslant (2k-3)(s-2)+1$. When $s \geqslant 2$ and $k \geqslant 3$, we have $|V_1| \geqslant 1$.

Step (2): Since $|V_1| \geqslant 1$, we can choose one vertex from $V_1$ and let this vertex be $w_2$. Let $E_2 = E_1 \cup \{l_{xw_2}, l_{yw_2}\}$ and $V_2 = V_1 \backslash \{l_{xw_2}, l_{yw_2}\}$. Similarly, we have $|V_2| \geqslant |V_1| - (2k-3) \geqslant (2k-3)(s-3)+1$.

Step ($r$): Let $3 \leqslant r \leqslant s-1$. Since $|(l_{xw_{r-1}} \cup l_{yw_{r-1}}) \cap V_{r-2}| \leqslant 2k-3$, we obtain $|V_{r-1}| = |V_{r-2} \backslash \{l_{xw_{r-1}}, l_{yw_{r-1}}\}| \geqslant |V_{r-2}| - (2k-3) \geqslant (2k-3)(s-r+1)+1 - (2k-3) = (2k-3)(s-r)+1 \geqslant 1$. Thus, from $V_{r-1}$, we can choose one vertex and denote it by $w_r$. Let $E_r = E_{r-1} \cup \{l_{xw_r}, l_{yw_r}\}$ and $V_r = V_{r-1} \backslash \{l_{xw_r}, l_{yw_r}\}$. Similarly, we get $|V_r| \geqslant |V_{r-1}| - (2k-3) \geqslant (2k-3)(s-r-1)+1$.

When $r = s-1$, we have $|V_{s-1}| \geqslant 1$. Thus, from $V_{s-1}$, we can choose a vertex and denote it by $w_s$. Since $w_s$ is adjacent to $x$ and $y$, two edges (denoted by $l_{xw_s}$ and $l_{yw_s}$) exist in $H$. Let $E_s = E_{s-1} \cup \{l_{xw_s}, l_{yw_s}\}$. Therefore, in $H$, $E_s$ contains $2s+1$ edges and these $2s+1$ edges form a Berge-$B_s$ in $H$. This contradicts the fact that $H$ is Berge-$B_s$-free. Thus, for any two adjacent vertices $x, y \in V(H)$, we obtain $|N_{xy}| \leqslant (k-1)(2s-1)-s$. $\square$

**Claim 4.2.** For any vertex $y \in V(H)$ and any vertex $x \in N_y$, we get $\sum_{i=2}^{k} e_i^x(N_y) \leqslant 2(k-1)(s-1)$, where $k \geqslant 3$ and $s \geqslant 2$.

**Proof.** For any vertex $y \in V(H)$ and any vertex $x \in N_y$, we will firstly prove $\sum_{i=2}^{k} e_i^x(N_y) \leqslant (2k-3)(s-1)+1$ by contradiction. Otherwise, we suppose that there exists a vertex $x \in N_y$ such that $\sum_{i=2}^{k} e_i^x(N_y) \geqslant (2k-3)(s-1)+2$. This means $|\bigcup_{i=2}^{k} E_i^x(N_y)| \geqslant (2k-3)(s-1)+2$. Since $x \in N_y$, we get $l_{xy} \in E(H)$. Obviously, $l_{xy} \in E_{k-1}^x(N_y)$. Thus, we get $|\left(\bigcup_{i=2}^{k} E_i^x(N_y)\right) \backslash l_{xy}| \geqslant (2k-3)(s-1)+1$. Namely, except for $l_{xy}$, $H$ has at least $(2k-3)(s-1)+1$ edges, each of which contains $x$ and at least two vertices in $N_y$.

Let $E^\star = \{f_i \mid f_i \in E(H), |f_i \cap N_y| \geqslant 2, x \in f_i, f_i \neq l_{xy}, 1 \leqslant i \leqslant (2k-3)(s-1)+1\}$. For any edge $f_i \in E^\star$, since $|f_i \cap N_y| \geqslant 2$ and $x \in f_i$, among $f_i$, we can choose a vertex (denoted by $z_i$) which belongs to $N_{xy}$, where $z_i \neq x, y$ and $1 \leqslant i \leqslant (2k-3)(s-1)+1$. Since $H$ is a linear hypergraph, $z_1, z_2, \cdots, z_{(2k-3)(s-1)+1}$ are different from each other.

Since $z_i \in N_{xy}$ and $l_{xy} \backslash \{x, y\} \subseteq N_{xy}$, we get $|N_{xy}| \geqslant (2k-3)(s-1) + 1 + k - 2 = (k-1)(2s-1) - s + 1$. This contradicts Claim 4.1. Therefore, for any vertex $y \in V(H)$ and any vertex $x \in N_y$, we have $\sum_{i=2}^{k} e_i^x(N_y) \leqslant (2k-3)(s-1) + 1 \leqslant 2(k-1)(s-1)$, where $s \geqslant 2$. $\qquad\square$

Let $B_u = V(H) \backslash (N_u \cup \{u\})$. For $v \in N_u$, let $S_v = \{w \mid w \in B_u, \exists l \in E_1^v(N_u), w \in l\}$. By the definition of $S_v$, we get $S_v \subseteq B_u$. We denote, by $l_i = \{u, i_2, \cdots, i_k\}$, the $i$-th edge which is incident with $u$. Furthermore, let $V_i = \bigcup_{j=2}^{k} S_{i_j}$, where $i = 1, 2, \cdots, d_u$.

**Claim 4.3.** $|V_i| \geqslant \sum_{j=2}^{k} |S_{i_j}| - \binom{k-1}{2}(2k-3)(s-1)$, where $1 \leqslant i \leqslant d_u$.

**Proof.** For any $i$ with $1 \leqslant i \leqslant d_u$ and $2 \leqslant p < q \leqslant k$, it follows from Claim 4.1 that $|S_{i_p} \cap S_{i_q}| \leqslant |N_{i_p i_q} \backslash (N_u \cup \{u\})| \leqslant |N_{i_p i_q} \backslash (l_i \setminus \{i_p, i_q\})| \leqslant (k-1)(2s-1) - s - (k-2) = (2k-3)(s-1)$. Thus, we get

$$
\begin{aligned}
|V_i| = |\bigcup_{j=2}^{k} S_{i_j}| &\geqslant \sum_{j=2}^{k} |S_{i_j}| - \sum_{2 \leqslant p < q \leqslant k} |S_{i_p} \cap S_{i_q}| \\
&\geqslant \sum_{j=2}^{k} |S_{i_j}| - \binom{k-1}{2}(2k-3)(s-1). \qquad\square
\end{aligned}
$$

**Claim 4.4.** We have $\sum_{i=1}^{d_u} |V_i| \leqslant (t-1)\big[n - (k-1)d_u - 1\big]$, where $t \geqslant 1$.

**Proof.** For any $v \in N_u$, since $S_v \subseteq B_u$, we have $V_i = \bigcup_{j=2}^{k} S_{i_j} \subseteq B_u$, where $1 \leqslant i \leqslant d_u$. Thus, we get $\bigcup_{i=1}^{d_u} V_i \subseteq B_u$. Next, we consider $\sum_{i=1}^{d_u} |V_i|$. For any vertex in $B_u$, we will prove that it contributes at most $t - 1$ times in $\sum_{i=1}^{d_u} |V_i|$. Otherwise, we suppose that there exists a vertex $w \in B_u$ such that it contributes at least $t$ times in $\sum_{i=1}^{d_u} |V_i|$. This means that at least $t$ sets of $V_1, V_2, \cdots, V_{d_u}$ contain $w$ simultaneously.

Without loss of generality, we suppose $w \in V_i$ with $i = 1, 2, \cdots, t$. According to the definition of $V_i$, there exists an edge $f_i \in E_1(N_u)$ such that $f_i \cap (l_i \backslash u) \neq \phi$ and $w \in f_i$, where $1 \leqslant i \leqslant t$. Let $1 \leqslant i < j \leqslant t$. We will prove $f_i \neq f_j$ for any $1 \leqslant i < j \leqslant t$. Otherwise, if there exist $1 \leqslant i < j \leqslant t$ such that $f_i = f_j$, then $f_i \cap (l_i \backslash u) \neq \phi$ and $f_i \cap (l_j \backslash u) \neq \phi$. Since $H$ is a linear hypergraph and $l_i \cap l_j = \{u\}$, we have $(l_i \backslash u) \cap (l_j \backslash u) = \phi$. Furthermore, since $l_i \backslash u, \ l_j \backslash u \subseteq N_u$, we obtain $|f_i \cap N_u| \geqslant 2$. This contradicts $f_i \in E_1(N_u)$, where $1 \leqslant i \leqslant t$. Therefore, $f_1, \ f_2, \cdots, f_t$ are different from each other. Thus, $f_1, \ f_2, \cdots, f_t, \ l_1, \ l_2, \cdots, l_t$ form a Berge–$K_{2,t}$ in $H$. This contradicts that $H$ is Berge–$K_{2,t}$–free. Finally, we get that each vertex in $B_u$ contributes at most $t - 1$ in $\sum_{i=1}^{d_u} |V_i|$. Furthermore, since $|B_u| =$ $n - (k-1)d_u - 1$, we have $\sum_{i=1}^{d_u} |V_i| \leqslant (t-1)(n-(k-1)d_u-1)$. Thus, we get Claim 4.4. $\square$

From Claim 4.2, for any $v \in N_u$, we have $d_v = \sum_{i=1}^{k} e_i^v(N_u) \leqslant e_1^v(N_u) + 2(k-1)(s-1)$. Therefore, for any $v \in N_u$, $e_1^v(N_u) \geqslant d_v - 2(k-1)(s-1)$. Since $H$ is a connected linear $k$-graph, we have $|S_v| = (k-1)e_1^v(N_u)$. Thus, $|S_v| \geqslant (k-1)d_v - 2(k-1)^2(s-1)$. Bearing $|N_u| = (k-1)d_u$ in mind, we get

$$
\begin{aligned}
\sum_{v \in N_u} |S_v| &\geqslant \sum_{v \in N_u} \left[ (k-1)d_v - 2(k-1)^2(s-1) \right] \\
&= \sum_{v \in N_u} (k-1)d_v - 2(k-1)^2(s-1)|N_u| \\
&= \sum_{v \in N_u} (k-1)d_v - 2(k-1)^3(s-1)d_u. \tag{13}
\end{aligned}
$$

By (13) and the fact $\sum_{v \in N_u} |S_v| = \sum_{i=1}^{d_u} \sum_{j=2}^{k} |S_{i_j}|$, we have

$$
\sum_{v \in N_u} (k-1)d_v \leqslant \sum_{i=1}^{d_u} \sum_{j=2}^{k} |S_{i_j}| + 2(k-1)^3(s-1)d_u. \tag{14}
$$

Furthermore, from Claims 4.3 and 4.4, (14) can be transformed into

$$
\begin{aligned}
\sum_{v \in N_u} (k-1)d_v &\leqslant \sum_{i=1}^{d_u} \left[ |V_i| + \binom{k-1}{2}(2k-3)(s-1) \right] + 2(k-1)^3(s-1)d_u \\
&= \sum_{i=1}^{d_u} |V_i| + \binom{k-1}{2}(2k-3)(s-1)d_u + 2(k-1)^3(s-1)d_u \\
&\leqslant (t-1)(n-1) + \Big[ \frac{1}{2}(k-1)(k-2)(2k-3)(s-1) \\
&\quad + 2(k-1)^3(s-1) - (k-1)(t-1) \Big] d_u.
\end{aligned}
$$

Therefore, we have

$$
\sum_{v \in N_u} d_v \leqslant \frac{(t-1)(n-1)}{k-1} + \left[ \frac{1}{2}(6k^2 - 15k + 10)(s-1) - (t-1) \right] d_u.
$$

By Corollary 2.2, we get

$$
\begin{aligned}
\rho_\alpha^2(H) &\leqslant \alpha d_u^2 + \frac{1-\alpha}{k-1} \sum_{v \in N_u} d_v \\
&\leqslant \alpha d_u^2 + (1-\alpha) \left[ \frac{(t-1)(n-1)}{(k-1)^2} + \frac{(6k^2 - 15k + 10)(s-1) - 2(t-1)}{2(k-1)} d_u \right].
\end{aligned}
$$

Therefore, we obtain Theorem 4.1. $\qquad\qquad\qquad\qquad\qquad\qquad\qquad\qquad\square$

It is noted that when $(6k^2 - 15k + 10)(s - 1) - 2(t - 1) \geqslant 0$, namely, $t \leqslant \frac{1}{2}(6k^2 - 15k + 10)(s - 1) + 1$, the value of the right-hand side of (12) is related with $d_u$. Therefore, we obtain Theorem 4.2 as follows.

**Theorem 4.2** *Let $H$ be a connected linear $k$-graph on $n$ vertices with the maximum degree $\Delta$, where $n, k \geqslant 3$. If $H$ is Berge-$\{B_s, K_{2,t}\}$–free, where $s \geqslant 2$ and $1 \leqslant t \leqslant \frac{1}{2}(6k^2 - 15k + 10)(s - 1) + 1$, then*

$$\rho_\alpha^2(H) \leqslant \alpha\Delta^2 + (1 - \alpha)\left[\frac{(t - 1)(n - 1)}{(k - 1)^2} + \frac{(6k^2 - 15k + 10)(s - 1) - 2(t - 1)}{2(k - 1)}\Delta\right].$$

If $H$ is a connected linear $k$-graph on $n$ vertices, then its maximum degree is less than or equal to $\frac{n-1}{k-1}$, where $n, k \geqslant 3$. Therefore, by Theorem 4.2, for $H \in \mathcal{H}_{n,k,s,t}$ with $n, k \geqslant 3$, $s \geqslant 2$ and $1 \leqslant t \leqslant \frac{1}{2}(6k^2 - 15k + 10)(s - 1) + 1$, we get the upper bounds of $\rho(H)$ and $\rho_\alpha(H)$, which are shown in Theorems 4.3 and 4.4, respectively.

**Theorem 4.3** *Let $H$ be a connected linear $k$-graph on $n$ vertices, where $n, k \geqslant 3$. If $H$ is Berge-$\{B_s, K_{2,t}\}$–free, where $s \geqslant 2$ and $1 \leqslant t \leqslant \frac{1}{2}(6k^2 - 15k + 10)(s - 1) + 1$, then*

$$\rho(H) \leqslant \sqrt{\frac{(6k^2 - 15k + 10)(s - 1)(n - 1)}{2(k - 1)^2}}.$$

**Theorem 4.4** *Let $H$ be a connected linear $k$-graph on $n$ vertices, where $n, k \geqslant 3$. If $H$ is Berge-$\{B_s, K_{2,t}\}$–free, where $s \geqslant 2$ and $1 \leqslant t \leqslant \frac{1}{2}(6k^2 - 15k + 10)(s - 1) + 1$, then*

$$\rho_\alpha^2(H) \leqslant \alpha\frac{(n - 1)^2}{(k - 1)^2} + (1 - \alpha)\frac{(6k^2 - 15k + 10)(s - 1)(n - 1)}{2(k - 1)^2}.$$

### Acknowledgments

The work was supported by the Natural Science Foundation of Shanghai under the grant number 21ZR1423500.

---